\documentclass[12pt]{iopart}

\usepackage{iopams} 
\usepackage{graphicx}
\usepackage{times}
\usepackage{epsfig}
\expandafter\let\csname equation*\endcsname\relax
\expandafter\let\csname endequation*\endcsname\relax
\usepackage{amsmath}
\usepackage{amssymb}
\usepackage{verbatim}
\usepackage{multirow}
\usepackage{makecell}
\usepackage{wrapfig} 
\usepackage{xcolor}
\usepackage{subfig}
\usepackage{epstopdf}
\usepackage{booktabs}
\usepackage{graphicx}
\usepackage{pgfplots}
\pgfplotsset{compat=1.6}
\usepackage{bigints} 
\usepackage[sort]{natbib}
\usepackage{float}
\usepackage[normalem]{ulem}
\usepackage{mathtools}

\begin{document}
\title[A Generalized Framework for Analytic Regularization]{A Generalized Framework for Analytic Regularization of Uniform Cubic B-spline Displacement Fields}

\author{Keyur D. Shah$^1$, James A. Shackleford$^1$, Nagarajan Kandasamy$^1$, Gregory C. Sharp$^2$}

\address{$^1$ Electrical and Computer Engineering Department, Drexel University, Philadelphia, PA 19104, USA}
\address{$^2$ Department of Radiation Oncology, Massachusetts General Hospital, Boston, MA 02114, USA}
\ead{gcsharp@partners.org}
\vspace{10pt}

\begin{abstract}
Image registration is an inherently ill-posed problem that lacks the constraints needed for a unique mapping between voxels of the two images being registered. As such, one must regularize the registration to achieve physically meaningful transforms. The regularization penalty is usually a function of derivatives of the displacement-vector field, and can be calculated either analytically or numerically. The numerical approach, however, is computationally expensive depending on the image size, and therefore a computationally efficient analytical framework has been developed. Using cubic B-splines as the registration transform, we develop a generalized mathematical framework that supports five distinct regularizers: diffusion, curvature, linear elastic, third-order, and total displacement. We validate our approach by comparing each with its numerical counterpart in terms of accuracy. We also provide benchmarking results showing that the analytic solutions run significantly faster --- up to two orders of magnitude --- than finite differencing based numerical implementations.
\end{abstract}

\section{Introduction}
The goal of image registration is to find a geometric transformation between corresponding image data that brings them into a common coordinate frame. By fusing multiple images, a physician gains a more complete understanding of patient anatomy.  The images can be acquired using similar or different imaging modalities --- for example, CT or MRI --- and may represent different stages of growth or disease. A registration is called \emph{rigid} if the motion or change is limited to global rotations and translations, and \emph{deformable} when the registration includes complex local variations. Deformable registration is preferred over rigid when locally precise alignment is needed; for example, in image-guided surgery~\citep{hartkens2003measurement} and image-guided radiotherapy~\citep{zhang2007automatic}. 

\begin{figure}[!ht]
    \begin{center}
        \includegraphics[scale=0.4]{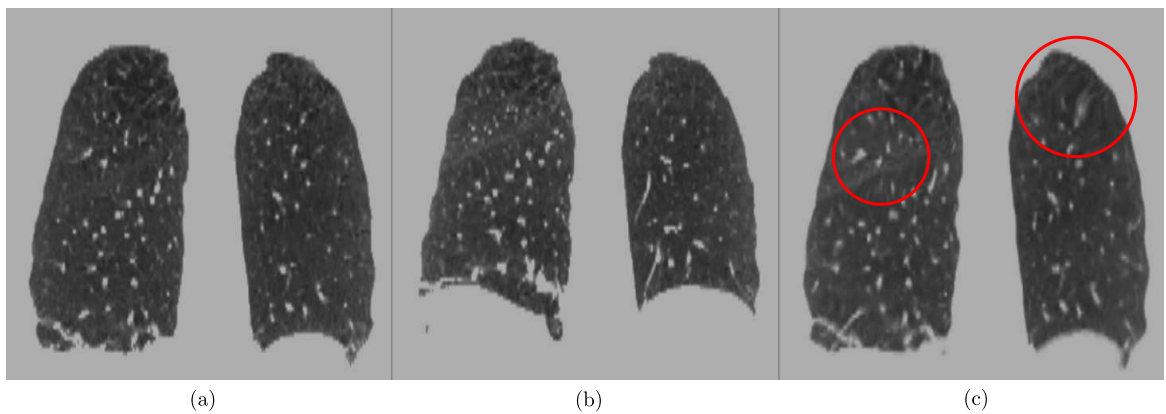}
        \caption{Example of intra-subject registration demonstrating the need for regularization.}
        \vspace{-12pt}
        \label{fig:intra}
    \end{center}
\end{figure}

Given a three-dimensional fixed image $F$ with voxel coordinates $\boldsymbol{x}=(x_1, x_2, x_3)$ and voxel intensity $f = F(\boldsymbol{x})$, and moving image $M$ with voxel coordinates $\boldsymbol{x'} = (x_1', x_2', x_3')$ and voxel intensity $m = M(\boldsymbol{x'})$ representing the same underlying anatomy as $F$ within the image overlap domain $\boldsymbol{\Omega}$, the two images $F$ and $M$ are said to be registered when the cost function 
\vspace*{-4pt}
\begin{equation}
C = \sum_{\boldsymbol{T}(\boldsymbol{x})\in\boldsymbol{\Omega}}{\boldsymbol{\Psi}\left(f,m\right)}
\label{eqn:general-cost}
\end{equation}
\vspace*{-4pt}
\noindent is minimized according to a similarity metric $\boldsymbol{\Psi}$ under the coordinate mapping
$\boldsymbol{T}(\boldsymbol{x})=\boldsymbol{x}+\boldsymbol{\nu}$. 
Here $\boldsymbol{\nu}$ is the
dense displacement field defined for every voxel $\boldsymbol{x}\in\boldsymbol{\Omega}$,
which is assumed capable of providing a good diffeomorphism from $F$ to $M$. A diffeomorphism is a globally one-to-one smooth and continuous mapping with derivatives that are invertible~\citep{ashburner2007fast}. Deformable image registration is an inherently ill-posed problem and so the unconstrained formulation in \eqref{eqn:general-cost} can lead to physically unrealistic transforms such as the one shown in Figure~\ref{fig:intra}. Moving image, Figure~\ref{fig:intra} (b) is registered to the fixed image,  Figure~\ref{fig:intra} (a). However, the resulting warped image Figure~\ref{fig:intra} (c) exhibits areas of irregular compression and expansion that are marked using red circles. 

It is desirable to confine the solution space to prevent physically unrealistic transforms. This can be done by introducing a penalty term which regularizes the transformation. We can modify \eqref{eqn:general-cost} to include the regularization term as 	

\begin{equation}
C=\sum_{\boldsymbol{T}(\boldsymbol{x})\in\boldsymbol{\Omega}}{\boldsymbol{\Psi}\left(f,m\right)} + S(\boldsymbol{\nu}),
\label{eqn:general-cost_modified}
\end{equation}
\noindent where the smoothness $S(\boldsymbol{\nu})$ is added to $\boldsymbol{\Psi}$ to drive $\boldsymbol{T}$ to a physically meaningful coordinate mapping. Several formulations for the regularization term $S$ have been developed in the literature.~\cite{rueckert1999nonrigid} penalize high bending energy in thin-plate spines to achieve smoother local deformations. This technique is implemented within the Medical Image Registration Toolkit (MIRTK).~\cite{Rohlfing2003} penalize local deviations from a unity Jacobian determinate to preserve incompressibility of volume regions since soft tissue is incompressible for small deformations. Linear elastic energy is minimized by~\cite{miller1993mathematical} to ensure that the deformation field generated by the registration is physically smooth.~\cite{chun2009simple} develop a regularizer based on sufficient conditions to enforce local invertibility.~\cite{cahill2009demons} generalize the Demon's algorithm~\citep{thirion1998image} to allow image-driven locally adaptive regularization.~\cite{sorzano2005elastic} propose a regularizer based on the curl and divergence of the underlying displacement field as it is the measure of the true roughness of the deformation. Use of Fourier methods to solve partial differential equations of various standard regularizers is studied by~\cite{cahill2007fourier}.~\cite{burger2013hyperelastic} develop a regularizer based on hyper-elasticity in the context of a mass-preserving registration problem. Using a Demons' framework,~\cite{tustison2013explicit} use the directly manipulated free-form deformation as a regularizer for the resulting displacement field to provide biologically plausible solutions.~\cite{vishnevskiy2016isotropic} use an isotropic version of the total variation regularization to correctly represent non-smooth displacement fields, that occur at sliding interfaces in the thorax and abdomen in image time-series during respiration.~\cite{schmidt2012estimation} and~\cite{delmon2013registration} proposed a direction-dependent regularization to estimate slipping organ motion.~\cite{miura2017impact} use the anatomically constrained deformation algorithm (ANACONDA) of RayStation, which penalizes invertibility of the displacement field and any large shape deviations to the region of interest.~\cite{fu2018adaptive} develop an adaptive direction dependent regularization technique using a Gaussian isotropic filter and a bilateral filter in order to preserve sliding motion.~\cite{ghaffari2017image} proposed a rank-regularized sum of squared differences similarity measure in order to overcome the challenge of spatially varying intensity distortion.~\cite{mang2016constrained} constrain the divergence of the velocity field to control the compressibility of the displacement field for a 2D case. 

Returning to the regularized cost function in~\eqref{eqn:general-cost_modified}, the penalty term {$S(\boldsymbol{\nu})$} can be calculated either numerically or analytically. The numerical approach, while simple and flexible, is computationally expensive depending on the image size, whereas the analytical methods requires solving a series of logical steps, but is computationaly inexpensive to solve. Moreover, the numerical methods provide the approximate solution and the analytical method gives the exact solution  and therefore analytic methods are often preferred.~\cite{shackleford2012analytic} have previously developed an analytic method to calculate the bending energy of a vector field that is parameterized via uniform cubic B-spline basis function and report significant speedup compared to the numerical counterpart. Along similar lines,~\cite{shusharina2012analytic} present an analytic method to regularize the radial basis function. 

While the above prior work has developed analytic methods for specific types of regularizers, the novelty of this paper lies in the development of a \emph{generalized mathematical framework} for this problem. Using cubic B-splines as the deformable registration transform, our framework accommodates five distinct types of regularization: diffusion~\citep{thirion1998image}, curvature~\citep{modersitzki2004numerical}, linear elastic~\citep{broit1981optimal}, third-order~\citep{lellmann2013anisotropic}, and total displacement. In addition to exhibiting improved computational efficiency with respect to numerical approaches, a key advantage of our approach is the ability to seamlessly combine multiple regularizers to realize custom or domain-specific smoothness constraints as dictated by the underlying registration problems at no additional cost in performance. 

The paper is organized as follows. Section~\ref{sec:background} discusses relevant theory and the regularizers of interest. Section~\ref{sec:derivation} develops the framework and its analytical solution; which is validated and benchmarked in Section~\ref{sec:results}. We conclude the paper in Section~\ref{sec:conclusions}.

\section{Background}\label{sec:background}
Here we describe the formulation of the diffusion, curvature, linear elastic, third-order, and total displacement regularizers for a three-dimensional displacement field parameterized by a uniform cubic B-spline basis function.

For volumetric or 3D registration, the displacement field at any given voxel is determined by the $4^3 = 64$ control points in the immediate vicinity of the voxel.  We use the term \emph{tile} to denote the set of voxels which receives local support from the same set of 64 control points. The tile forms the backbone of an analytic expression for the continuous displacement field  $\boldsymbol{\nu}$. B-spline interpolation is performed for each vector within a tile using the 64 control-point coefficients that provide local support for the operation. The B-spline interpolation yielding the first component of the displacement vector for a voxel located at $\boldsymbol{x}$ is
\begin{equation}
\nu_{1}(\boldsymbol{x}) = \sum_{l=0}^3 \sum_{m=0}^3 \sum_{n=0}^3 \beta_{l}(u_1)\beta_{m}(u_2)\beta_{n}(u_3) p_{1,l,m,n}
\label{Displacement vector field}
\end{equation}
\noindent where $p_{1}$ is {one of the 64} B-spline coefficient used to interpolate the $\nu_1$ component of the displacement vector $\boldsymbol{\nu}$ for the voxel located at $\boldsymbol{x}$. The $\beta_l$, $\beta_m$, and $\beta_n$ terms represent the uniform cubic B-spline basis functions in the $x_1$, $x_2$, and $x_3$ directions, respectively, and $u_1$, $u_2$, and $u_3$ represent the normalized local coordinates of the voxel within its tile~\citep{shackleford2010developing}. The uniform cubic B-spline basis function $\beta_l$ along the $x_1$ direction is given by
\begin{align}
\beta_l(u_1) =
\begin{dcases}
\frac{(1 - u_1)^3}{6} & l = 0 \\
\frac{3{u_1}^3 - 6{u_1}^2 + 4}{6} & l = 1\\
\frac{-3{u_1}^3 + 3{u_1}^2 + 3u_1 + 1}{6} & l = 2\\ \label{eqn:cubic}
\frac{{u_1}^3}{6} & l = 3
\end{dcases}
\end{align}
and similarly for $\beta_m$ and $\beta_n$ in the $x_2$ and $x_3$ directions respectively.
The displacement vector components $\nu_2$ and $\nu_3$ in the remaining two directions are calculated similarly. The optimizer updates the B-spline coefficients during each iteration until an optimal registration between the moving and fixed images is achieved. The regularization penalty term is a function of the displacement vector field \emph{i.e.} a function of the B-spline coefficients. Thus, the regularization penalty term for the entire deformation may be expressed as a sum of the regularization penalty terms over all tiles.  Therefore, our approach is to develop an operator that computes the penalty term for a tile as a function of its B-spline control points.

We now describe the five regularizers of interest. The functions in~\eqref{eqn:functions1} describe the first, second, and third-order partial derivatives of the displacement field, which are the building blocks of these regularizers:
\begin{align}
f_1(\nu_i,x_j) &= \frac{\partial{\nu_{i}}}{\partial{x}_j} \nonumber \\
f_2(\nu_i,x_j,x_k) &= \frac{\partial^2\nu_i}{\partial{x}_j \partial{x}_k} \nonumber \\
f_3(\nu_i,x_j,x_k) &= \frac{\partial{\nu_{i}}}{\partial{x}_j} \frac{\partial{\nu_{i}}}{\partial{x}_k} \nonumber \\ 
f_4(\nu_i,x_j,x_k,x_q) &= \frac{\partial^3\nu_i}{\partial{x}_j \partial{x}_k \partial{x}_q} \label{eqn:functions1}
\end{align}

\vspace{6pt}

\noindent \textbf{Diffusion Regularizer:} Thirion originally  introduced the demons algorithm as a diffusion process~\citep{thirion1998image} and in later work Modersitzki coined the term diffusion regularizer~\citep{modersitzki2004numerical} since this gradient-based partial differential equation was viewed as a generalized diffusion equation. The penalty term is given by
\begin{equation}
\begin{aligned}
S_{1} = \int_{\boldsymbol{\Omega}} { \sum_{i,j=1}^3f_1(\nu_i,x_j)^2}d\boldsymbol{x},
\label{eqn:smooth_diff}
\end{aligned}
\end{equation}

\noindent where $\boldsymbol{\Omega}$ denotes the image region over which the regularization penalty is to be calculated.

\noindent \textbf{Curvature Regularizer:} A curvature regularizer uses second-order derivative terms and aims to reduce the bending energy of the displacement field~\citep{modersitzki2004numerical}. The penalty term is specified as 
\begin{equation}
\begin{aligned}
S_{2} = \int_{\boldsymbol{\Omega}} \sum_{i,j,k=1}^3f_2(\nu_i,x_j,x_k)^2d\boldsymbol{x}.
\label{eqn:smooth_curv}
\end{aligned}
\end{equation}

\noindent \textbf{Linear Elastic Regularizer:} The penalty term is given by
\begin{equation}
\begin{aligned}
S_{3} = \int_{\boldsymbol{\Omega}} { \sum_{i,j=1}^3f_1(\nu_i,x_j)^2} + { \sum_{i,j,k=1, j \neq k}^3f_3(\nu_i,x_j,x_k)}d\boldsymbol{x},
\label{eqn:smooth_elastic}
\end{aligned}
\end{equation}
which aims to strike a balance between the global registration achieved via affine mapping versus the more local elastic registration~\citep{broit1981optimal}.

\vspace{6pt}

\noindent \textbf{Third-order Regularizer:} The penalty term is specified in terms of third-order derivative terms as
\begin{equation}
\begin{aligned}
S_{4} = \int_{\boldsymbol{\Omega}} \sum_{i,j,k,o=1}^3f_4(\nu_i,x_j,x_k,x_o)^2 d\boldsymbol{x}.
\label{eqn:smooth_to}
\end{aligned}
\end{equation}

\vspace{6pt}

\noindent \textbf{Total displacement Regularizer:} The penalty term is specified in terms of the magnitude of the displacement vector field at each voxel as
\begin{equation}
\begin{aligned}
S_{5} = \int_{\boldsymbol{\Omega}} \sum_{i=1}^3{\nu_{i}}^2 d\boldsymbol{x}.
\label{eqn:smooth_td}
\end{aligned}
\end{equation}

\noindent Combining the above-described regularizers, the total smoothness penalty $S$ for the unified framework can be expressed as
\begin{align}
S = \mu_1 S_{1} + \mu_2 S_{2} + \mu_3 S_{3} + \mu_4 S_{4} + \mu_5 S_{5},
\label{eqn:smooth}
\end{align}
where $\mu_1$, $\mu_2$, $\mu_3$, $\mu_4$, and $\mu_5$ are the weights corresponding to the diffusion, curvature, linear-elastic, third-order, and total displacement regularizers, respectively. In the simplest case, one regularizer can be chosen over the others; for example, setting $\mu_2$, $\mu_3$, $\mu_4$, and $\mu_5$ to zero chooses the diffusion regularizer. Other regularizers can be selected similarly. 
Several other regularizers can be derived using this framework as long as the regularizer penalty terms consists of product of two partial derivative terms of the displacement filed from ~\eqref{eqn:functions1} or one of the components of the displacement field.

\section{Development of Analytical Algorithm}\label{sec:derivation}
  
The analytic algorithm developed in this section comprises the following three major steps: (1) taking the first, second, or third-order derivative of the displacement vector field; (2) squaring or multiplying the derivative terms; and (3) integrating the products of the derivative terms over a tile. The various derivative terms that constitute \eqref{eqn:smooth} are recast as simple matrix operations that can be efficiently performed using Basic Linear Algebra Subprograms (BLAS). The algorithmic steps are described in greater detail below. 
 
When represented sparsely via the uniform cubic B-spline basis, the displacement field $\boldsymbol{\nu}$ is parameterized by the set of B-spline basis coefficients $\boldsymbol{p_1},\boldsymbol{p_2},\boldsymbol{p_3}$, where

\begin{equation}
    \boldsymbol{p_1} = 
    \begin{bmatrix}
    p_{1,0,0,0} \\
    \vdots \\
    p_{1,I,J,K} \\
    \end{bmatrix},
    \boldsymbol{p_2} = 
    \begin{bmatrix}
    p_{2,0,0,0} \\
    \vdots \\
    p_{2,I,J,K} \\
    \end{bmatrix}, 
    \boldsymbol{p_3} = 
    \begin{bmatrix}
    p_{3,0,0,0} \\
    \vdots \\
    p_{3,I,J,K} \\
    \end{bmatrix}
\end{equation}
\noindent are the control points that define the displacement field within a single tile. For cubic B-splines, $I=J=K=2$.  The tile has dimensions of $\boldsymbol{r}=r_1 \times r_2 \times r_3$ mm$^3$. To express the first component of the displacement field $\nu_1$ as a function of $\boldsymbol{p_1}$, we start with the matrix $\boldsymbol{B}$ containing the coefficients for the cubic B-spline basis function as described in~\eqref{eqn:cubic} and matrix $\boldsymbol{R_1}$ which controls for tile size as
\begin{align}
\boldsymbol{B} = \frac{1}{6}
\begin{bmatrix}
1 & -3 &  3 & -1 \\
4 &  0 & -6 &  3 \\
1 &  3 &  3 & -3 \\
0 &  0 &  0 &  1
\end{bmatrix}, \;  
\boldsymbol{R_1} =
\begin{bmatrix}
1 & 0             & 0               & 0   \\
0 & \frac{1}{r_1} & 0               & 0   \\
0 & 0             & \frac{1}{r_1^2} & 0   \\
0 & 0             & 0               & \frac{1}{r_1^3}
\end{bmatrix}.
\label{eqn: Initialize1}
\end{align} 

\noindent We also generate the matrix $\boldsymbol{\Delta}^{(\delta)}$ which is defined for $\delta \in$ [0,3] as
\begin{align}
\boldsymbol{\Delta}^{(0)} = 
\begin{bmatrix}
1 & 0 & 0 & 0   \\
0 & 1 & 0 & 0   \\
0 & 0 & 1 & 0   \\
0 & 0 & 0 & 1
\end{bmatrix},  \; 
\boldsymbol{\Delta}^{(1)} = 
\begin{bmatrix}
0 & 0 & 0 & 0   \\
1 & 0 & 0 & 0   \\
0 & 2 & 0 & 0   \\
0 & 0 & 3 & 0
\end{bmatrix},  \; \nonumber \\
\boldsymbol{\Delta}^{(2)} = 
\begin{bmatrix}
0 & 0 & 0 & 0   \\
0 & 0 & 0 & 0   \\
2 & 0 & 0 & 0   \\
0 & 6 & 0 & 0
\end{bmatrix},   \; 
\boldsymbol{\Delta}^{(3)} = 
\begin{bmatrix}
0 & 0 & 0 & 0   \\
0 & 0 & 0 & 0   \\
0 & 0 & 0 & 0   \\
6 & 0 & 0 & 0
\end{bmatrix}.   \;
\label{eqn:derivs-p2}
\end{align}

\noindent Matrices $\boldsymbol{Q_1}$, $\boldsymbol{Q_2}$, and $\boldsymbol{Q_3}$ can now be calculated as 
\begin{align}
\boldsymbol{Q_1}^{(\delta)} = \boldsymbol{B} \boldsymbol{R_1} \boldsymbol{\Delta}^{(\delta)}, \; 
\boldsymbol{Q_2}^{(\delta)} = \boldsymbol{B} \boldsymbol{R_2} \boldsymbol{\Delta}^{(\delta)}, \; \mathrm{and} \; 
\boldsymbol{Q_3}^{(\delta)} = \boldsymbol{B} \boldsymbol{R_3} \boldsymbol{\Delta}^{(\delta)}, 
\label{eqn:derivs-p1}
\end{align}
\noindent to provide a convenient method for obtaining the first, second and third-order derivatives, $\boldsymbol{\nu'}$, $\boldsymbol{\nu''}$ and $\boldsymbol{\nu'''}$, respectively, with respect to the Euclidean basis as required by the calculation of the smoothness penalty. For example, $\boldsymbol{Q_1}^{(0)}$ is used when $\boldsymbol{\nu}$ is needed as is, and $\boldsymbol{Q_1}^{(1)}$, $\boldsymbol{Q_1}^{(2)}$ and $\boldsymbol{Q_1}^{(3)}$ are used when the first, second and third derivatives, respectively, of $\boldsymbol{\nu}$ are needed in the calculations. These matrices also map the domain of the B-spline basis function to lie in the interval [0,1]. Now, the vector $\nu_1$ may now be expressed at a point $\boldsymbol{x}={x_1,x_2,x_3}$ using the 64 B-spline coefficients supporting $\boldsymbol{x}$ as the tensor product:
\begin{align}
\nu_1 = \sum_{i,j,k=0}^{3} {p}_{i,j,k}
\sum_{a=0}^3 \boldsymbol{Q_1}^{(0)} \boldsymbol{x_1}
\sum_{b=0}^3 \boldsymbol{Q_2}^{(0)} \boldsymbol{x_2}
\sum_{c=0}^3 \boldsymbol{Q_3}^{(0)} \boldsymbol{x_3}
\label{eqn:tensor}
\end{align}
where the four rows of $\boldsymbol{Q_1}^{(\delta)}$ can be separated into vectors as 

\begin{equation}
\boldsymbol{Q_1}^{(\delta)} =
\begin{bmatrix}
q_{1,0}^{\boldsymbol{T}} \\
q_{1,1}^{\boldsymbol{T}} \\
q_{1,2}^{\boldsymbol{T}} \\
q_{1,3}^{\boldsymbol{T}}
\end{bmatrix}^{(\delta)}
\label{eqn:qrows}
\end{equation}

\noindent We define 
\begin{equation}
\boldsymbol{x_1} = [1\quad x_1\quad x_1^2\quad x_1^3]^{\boldsymbol{T}}
\label{eqn:cartesian-x}
\end{equation}
 
Cartesian basis vectors $\boldsymbol{x_2}$ and $\boldsymbol{x_3}$ are defined similarly. Finally, vectors $\nu_2$ and $\nu_3$ are defined in a similar fashion to $\nu_1$ to complete the B-spline interpolation operation.

Now that we have the basic mathematical framework in place, the next steps of the algorithm are designed to calculate the various squared and product terms listed in \eqref{eqn:smooth} --- specifically, to square or to multiply the $\sum \boldsymbol{Q_1}^{(\delta)}
\sum  \boldsymbol{Q_2}^{(\delta)} \sum \boldsymbol{Q_3}^{(\delta)}$ sub-term in \eqref{eqn:tensor}. 

\noindent A $4 \times 4$ matrix is calculated as the outer product

\begin{equation}
\boldsymbol{\Xi_{1,a,b}}\ = q_{1,a}^{\delta_{i}} \otimes q_{1,b}^{\delta_{j}},
\label{eqn: Xi}    
\end{equation}
\noindent where $a, b\in \{0, 1, 2, 3\}$ and $i, j \in \{1, 2, 3\}$. Here, $\delta_{i}$ and $\delta_{j}$ are the same when taking the square of the derivative term whereas they are different when taking the product of two distinct terms --- which is the case for the linear elastic regularizer. The $\delta$ terms indicate the variable on which the partial derivative of the displacement field is obtained. For example, if $i$ is 1, the partial derivative is taken with respect to the first component, and so on. The matrices $\boldsymbol{\Xi_{2,a,b}}$ and $\boldsymbol{\Xi_{3,a,b}}$ for the second and third components can be calculated similarly. For the ease of readability, the $\delta$ terms are not included in the follow-up equations.

Taking the outer product of the Cartesian basis vector $\boldsymbol{x_1}$ as defined in \eqref{eqn:cartesian-x} with itself results in a $4 \times 4$ matrix  
\begin{equation}
\boldsymbol{X_1} = \boldsymbol{x_1} \otimes \boldsymbol{x_1} = \begin{bmatrix}
    1 & x_1 & x_1^2 & x_1^3   \\
    x_1 & x_1^2 & x_1^3 & x_1^4   \\
    x_1^2 & x_1^3 & x_1^4 & x_1^5   \\
    x_1^3 & x_1^4 & x_1^5 & x_1^6
    \end{bmatrix},
\label{eqn:X-cartesian}    
\end{equation}
\noindent and $\boldsymbol{X_2}$ and $\boldsymbol{X_3}$ can be obtained similarly. Taking the Hadamard product of \eqref{eqn: Xi} and \eqref{eqn:X-cartesian} results in $4 \times 4$ matrix
\begin{equation}
\psi_{1,a,b}\ = \boldsymbol{\Xi_{1,a,b}}\ \odot \boldsymbol{X_1}.
\label{eqn: psi_x}
\end{equation}
\noindent The elements of $\psi_{1,a,b}$ can be combined into a matrix $\boldsymbol{\Gamma_1}$, where each element $\boldsymbol{\Gamma_1}(a,b)$ is formed as
\begin{equation}
\boldsymbol{\Gamma_1}(a,b) = \sum_{i=1}^4 \sum_{j=1}^4 \psi_{1,a,b}(i,j).
\label{eqn: Gamma}
\end{equation}
\noindent Since there are 16 combinations for $a$ and $b$, $\boldsymbol{\Gamma_1}$ is a $4 \times 4$ matrix; $\boldsymbol{\Gamma_2}$ and $\boldsymbol{\Gamma_3}$ are obtained similarly, allowing for the desired composite matrix operator 
\begin{equation}
    \boldsymbol{\Gamma} = \boldsymbol{\Gamma_1} \otimes \boldsymbol{\Gamma_2}\otimes \boldsymbol{\Gamma_3},
\label{eqn: Gamma1}
\end{equation}
\noindent to calculate the smoothness metric over a tile, for the specified choice of $\delta$'s, as
\begin{equation}
    \int_{0,0,0}^{r_1,r_2,r_3}\sum_{}
                    \boldsymbol{p_{i}}^{\boldsymbol{T}} \boldsymbol{\Gamma} \boldsymbol{p_{j}}d\boldsymbol{x}. 
\label{eqn:s-final}
\end{equation}
\noindent Since B-spline coefficients are constants, we can rewrite the above expression as 
\begin{equation}
   \sum_{}
                    \boldsymbol{p_{i}}^{\boldsymbol{T}} \boldsymbol{V} \boldsymbol{p_{j}}, 
\label{eqn:s-final1}
\end{equation}
\noindent where $\boldsymbol{V} = \int_{0,0,0}^{r_1,r_2,r_3} \boldsymbol{\Gamma} d \boldsymbol{x}$ and $i,j \in {1,2,3}$. The smoothness penalty $S$ for the entire volume can be calculated as the sum of all its constituent tiles.

\vspace{6pt}
This analytic implementation is extremely memory efficient. The 32 $\boldsymbol{V}$ matrices each having dimensions of $64 \times 64$ corresponding to every term in the five regularizers are calculated beforehand and reused in the entire optimization process. During the optimization step only $\boldsymbol{p^T} \boldsymbol{V} \boldsymbol{p}$ is calculated for each tile according to~\eqref{eqn:s-final1}. The overall memory requirement is 512 kB ($32 \times \boldsymbol{V} = 32 \times 64 \times 64 \times 4\ kB)$.

\noindent \textbf{Illustrative Examples:} Given the above-described formalism, the generalized equation for squaring or multiplying the various derivative terms can be written in terms of the B-spline basis coefficients and the composite matrix operator $\boldsymbol{\Gamma}$ as
\begin{equation}
{\left(\frac{\partial^n\nu_i}{\partial{x_i}}\frac{\partial^n\nu_j}{\partial{x_j}}\right)} = \boldsymbol{p_i}^{\boldsymbol{T}}  \boldsymbol{\Gamma}^{n\times(\boldsymbol{\delta_{i}},\boldsymbol{\delta_{j}})} \boldsymbol{p_j},
\label{eqn: General_example}
\end{equation}
\noindent where $n$ determines the order of the partial derivative and $\boldsymbol{\delta_i} = (\delta_{i1}, \delta_{j1}, \delta_{k1})$ and $\boldsymbol{\delta_j} = (\delta_{i2}, \delta_{j2}, \delta_{k2})$. Individual terms of $\boldsymbol{\delta_i}$ depend on the variable on which the partial derivative of the displacement field is obtained. For example, $\delta_{i1} = 1$ if the partial derivative is taken with respect to the first component. For the squared terms $\boldsymbol{\delta_i}$ and $\boldsymbol{\delta_j}$ are the same. 
The terms corresponding to the diffusion, curvature, linear-elastic, third-order, and total displacement regularizers are described by the following equations using the appropriate $\delta$ values discussed earlier in~\eqref{eqn: Xi} and \eqref{eqn: General_example}.
\begin{align}
{\left(\frac{\partial\nu_2}{\partial{x_1}}\right)}^2 &= \boldsymbol{p_2}^{\boldsymbol{T}}  \boldsymbol{\Gamma}^{(100,100)}  \boldsymbol{p_2}  \nonumber \\
{\left(\frac{\partial^2\nu_1}{\partial{x_1}\partial{x_3}}\right)}^2 &= \boldsymbol{p_1}^{\boldsymbol{T}} \boldsymbol{\Gamma}^{(101,101)}  \boldsymbol{p_1} \nonumber \\
{\left(\frac{\partial\nu_1}{\partial{x_1}}\frac{\partial\nu_3}{\partial{x_3}}\right)} &= \boldsymbol{p_1}^{\boldsymbol{T}}  \boldsymbol{\Gamma}^{(100,001)}  \boldsymbol{p_3} \\
{\left(\frac{\partial^3\nu_3}{\partial{x_1}^3}\right)}^2 &= \boldsymbol{p_3}^{\boldsymbol{T}} \boldsymbol{\Gamma}^{(300,300)}  \boldsymbol{p_3} \nonumber \\ 
{\left(\nu_3\right)}^2 &= \boldsymbol{p_3}^{\boldsymbol{T}} \boldsymbol{\Gamma}^{(000,000)} \boldsymbol{p_3} \nonumber
\label{eqn:illustrative_example}
\end{align}
\noindent Using similar reasoning, the linear elastic regularization penalty would be computed as the sum of twelve vector-matrix-vector products of the form $\boldsymbol{p^T} \boldsymbol{V} \boldsymbol{p}$.    
\vspace{-6pt}
\section{Experiments and Results}\label{sec:results}
\begin{figure}[!ht]
\begin{center}
\includegraphics[width=\textwidth]{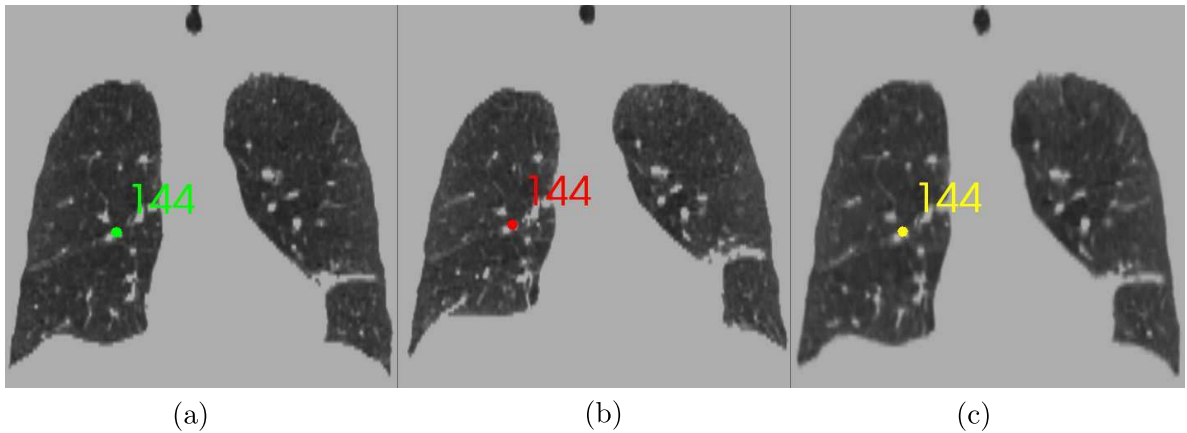}
\caption{Example of intra-subject registration showing a coronal slice of the (a) fixed, (b) moving and (c) warped images, with the corresponding selected landmark point overlaid.}
\label{fig:landmarks}
\end{center}
\vspace{-12pt}
\end{figure}

We quantify performance of the developed methods in terms of: accuracy and speedup. The DIR-Lab dataset used in our experiments consists of 4D-CT images from ten patients who were treated for malignancies in the esophagus or lung~\citep{castillo2009b}. The CT images were masked to only include the lungs, trachea and bronchi. The maximum inhalation phase (fixed image) was registered with the maximum exhalation phase (moving image) for intra-subject registrations. Five of the ten studies have volumes of $512 \times 512 \times 128$ voxels with physical separation of $0.92 \times 0.92 \times 2.5$ mm, and the remaining volumes have $256 \times 256 \times 94$ voxels. In each case, 300 landmark points were placed within the lung by a medical expert. Figure~\ref{fig:landmarks} shows a coronal slice for an intra-subject registration example with a landmark point overlaid. The moving landmarks are warped using the underlying transformation $\mathbf{T}$ to get the warped landmarks. Registration accuracy is measured as the average separation between the fixed and warped landmark points.

The regularizers developed in this paper have been implemented within Plastimatch, an open source software for image computation with focus on high-performance volumetric registration of medical images. Plastimatch is distributed under a BSD-style license and can be downloaded from {www.plastimatch.org}.
The fixed and moving images are registered using a three-stage pyramidal registration approach. Grid spacing for the first and second stage was kept constant at 100 mm and 80 mm, respectively. Grid spacing for the third stage was varied from a coarse value of 60 mm to a finer value of 10 mm. Registrations are performed using the mean-squared error similarity metric, penalized by $S_n$ with weight $\mu_n$. Referring to~\eqref{eqn:smooth}, we choose a single regularization strategy over others by setting the remaining weights to zero. For example, setting $\mu_1$, $\mu_3$ , $\mu_4$, and $\mu_5$ to zero chooses the curvature regularizer. The B-spline coefficients $\boldsymbol{P}$ describing the transform $\boldsymbol{T}$ are optimized via the L-BFGS-B optimizer using an analytically computed cost function and gradient. 

Accuracy is measured using the \emph{mean landmark separation} (\emph{MLS}) between the inhale landmarks and warped exhale landmarks (or vice versa) after registration, and smoothness is measured using the \emph{minimum Jacobian determinant} of the resulting displacement vector field over the entire volume. Experiments were performed to assess both quantities as a function of control-point spacing as well as $\mu_n$. 

\begin{figure*}[!ht]
    \centering 
    \subfloat[]{\label{fig:mls_masked}\includegraphics[width=\textwidth]{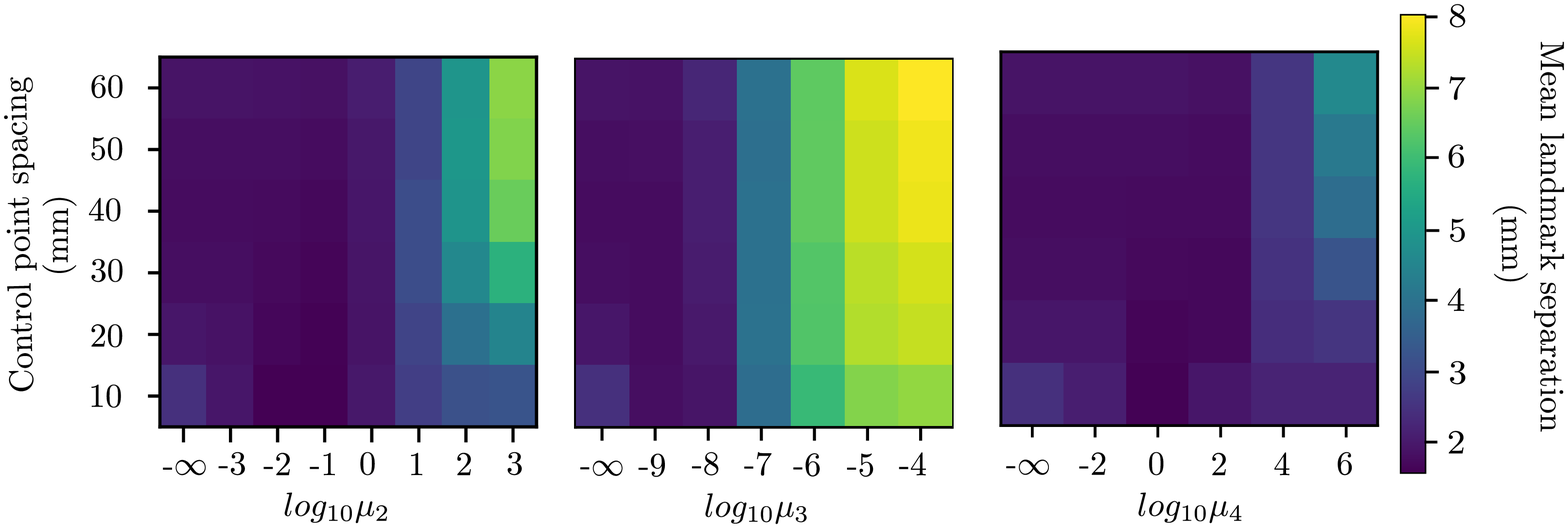}} \hfill     
    \subfloat[]{\label{fig:mjd_masked}\includegraphics[width=\textwidth]{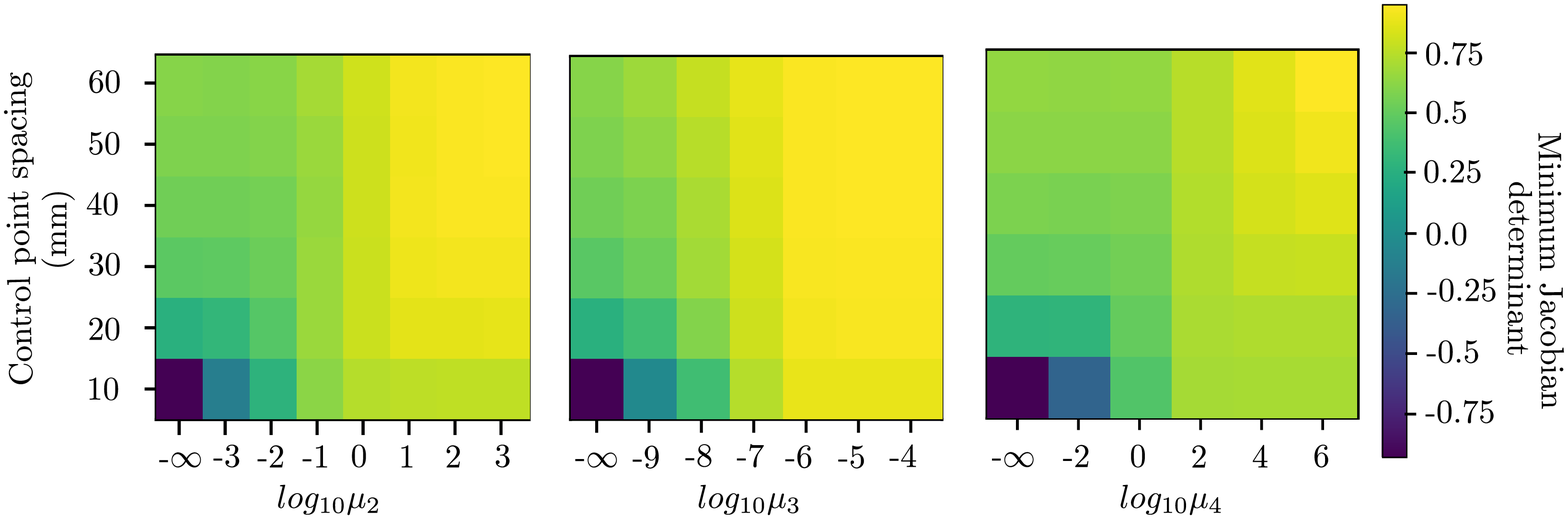}}
    \caption{Mean landmark separation (a) and and minimum Jacobian determinant of the transform $\boldsymbol{T}$ (b) as a function of B-spline control-point spacing and the weight $\mu_n$,  over the 10 thoracic cases using curvature (left), linear elastic (middle), and third-order (right) regularizer for the masked images.
    }
    \label{fig:graphs1}
\end{figure*}

Table~\ref{tab:my-table1} shows the relative difference (in \%) between the MLS achieved by the analytic and numeric implementations. The maximum relative difference of 7.4\% occurs for the third-order regularizer. This is because voxels located along the image boundaries are not used to calculate the smoothness penalty in case of the finite differencing numeric solutions.  Figure~\ref{fig:graphs1} shows accuracy and smoothness results for the curvature, linear elastic, and third-order regularizers. For a given control-point spacing, a desirable $\mu_n$ is one which produces the best compromise between a small MLS and smooth $\boldsymbol{T}$. A smooth $\boldsymbol{T}$ can be inferred by a positive minimum Jacobian determinant. A smooth nonrigid transformation ensures that the warped image after registration is free from unrealistic compression and/or expansion artifacts~\citep{chun2009simple}. Considering the curvature regularizer, for example, the MLS achieved when $\mu_2 = 10^{-3}$ is very similar without application of any regularization penalty, and smaller weights need not be explored. Conversely, the MLS achieved when  $\mu_2 = 10^3$ is nearly equal to the MLS prior to registration and changes little when $\mu_2$ is increased further. Thus, weights larger than $10^{3}$ are not shown. Lower and upper bounds for the weights of the remaining four regularizers are determined similarly. Additionally looking at the minimum Jacobian determinant heatmaps from Figure~\ref{fig:graphs1} (b), the increased need of regularization at smaller control point spacing can be inferred. The minimum Jacobian determinant of the resulting displacement vector field without any regularization penalty is negative for control point spacing of 10 mm for all the three regularizers shown. This is due to the higher degrees of freedom available to the displacement vector field at such smaller control point spacing.

\begin{table*}[!ht]
\caption{Relative difference (in \%) between MLS achieved by the proposed analytic schemes versus numerical finite differencing of the displacement field.}
\label{tab:my-table1}
\centering
\small
\begin{tabular}{@{}ccccccc@{}}
\toprule
Volume size                 & Grid spacing             & Diffusion      & \multicolumn{1}{c}{Curvature}      & \multicolumn{1}{c}{Elastic}  & \multicolumn{1}{c}{Third-order}     & Tot. displacement  \\ \midrule
$256 \times 256 \times 94$  & $20 \times 20 \times 20$ &     1.5     &             1.9                 &     1.3                   &              0.9                 &     4.4          \\
$256 \times 256 \times 94$  & $30 \times 30 \times 30$ &     0.1     &             3.1                 &     0.2                   &              3.9                 &     4.1          \\
$512 \times 512 \times 128$ & $20 \times 20 \times 20$ &     2.3     &             2.3                 &     3.0                   &              7.4                 &     1.9          \\
$512 \times 512 \times 128$ & $30 \times 30 \times 30$ &     5.9     &             2.0                 &     3.9                   &              5.7                 &     2.9          \\ \bottomrule
\end{tabular}
\end{table*}

\vspace{12pt}

\begin{table}[!ht]
\centering
\caption{Lowest MLS value achieved for each of the ten cases and the corresponding registration configuration.}
\label{tab:my-table3}
\begin{tabular}{@{}cccc@{}}
\toprule
Case \# & MLS (mm) & Regularizer ($\mu$)  & Grid Spacing (mm) \\ \midrule
1       & 1.15    &   Diffusion ($10^{-2}$)  & 30        \\
2       & 1.28    &   Diffusion ($10^{-3}$) & 10         \\
3       & 1.52    &   Linear Elastic ($10^{-8}$) & 10          \\ 
4       & 1.56    &   Curvature ($10^{-1}$) & 20          \\
5       & 1.80    &   Curvature ($10^{-2}$) & 10         \\
6       & 1.56    &   Curvature ($10^{-1}$) & 10        \\
7       & 1.74    &   Curvature ($10^{-2}$) & 10       \\
8       & 1.49    &   Curvature ($10^{-2}$) & 10       \\
9       & 1.50    &   Third-order ($1$)     & 10     \\
10      & 1.42    &   Curvature ($10^{-2}$) & 10       \\ \bottomrule
\end{tabular}
\end{table}
Table~\ref{tab:my-table3} summarizes the results for each of the ten cases by providing the least MLS achieved and the registration configuration used. Curvature regularizer works best for most of the cases, but other registration problems may be best optimized with a different regularization choice.
\begin{figure*}[!ht]
    \centering
    \subfloat[Image overlay]{\label{fig:difference}\includegraphics[height=1.5in,width=0.24\textwidth]{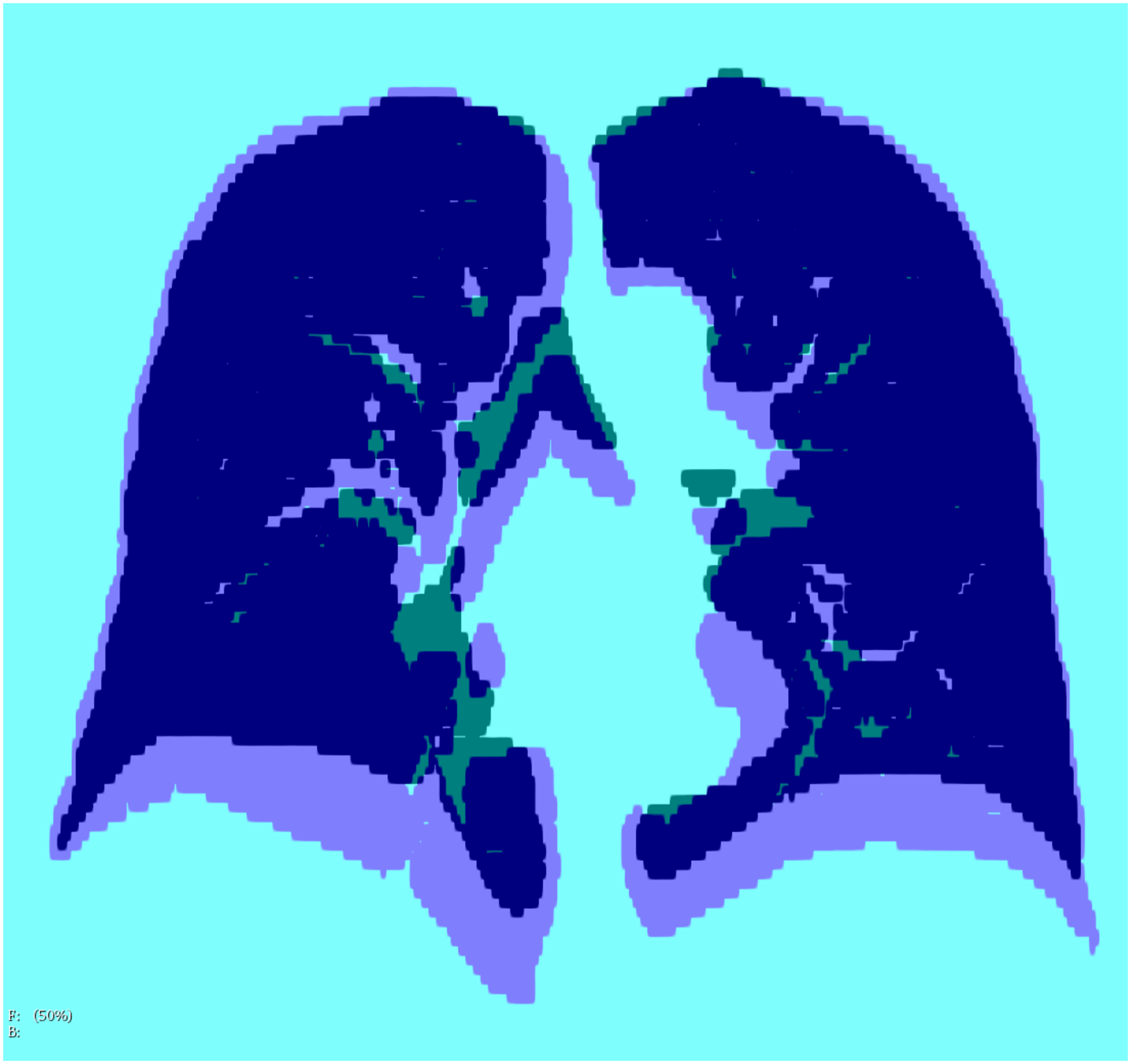}} \hfill
    \subfloat[$\boldsymbol{T}$ with no regularization]{\label{fig:no_reg}\includegraphics[height=1.5in,width=0.24\textwidth]{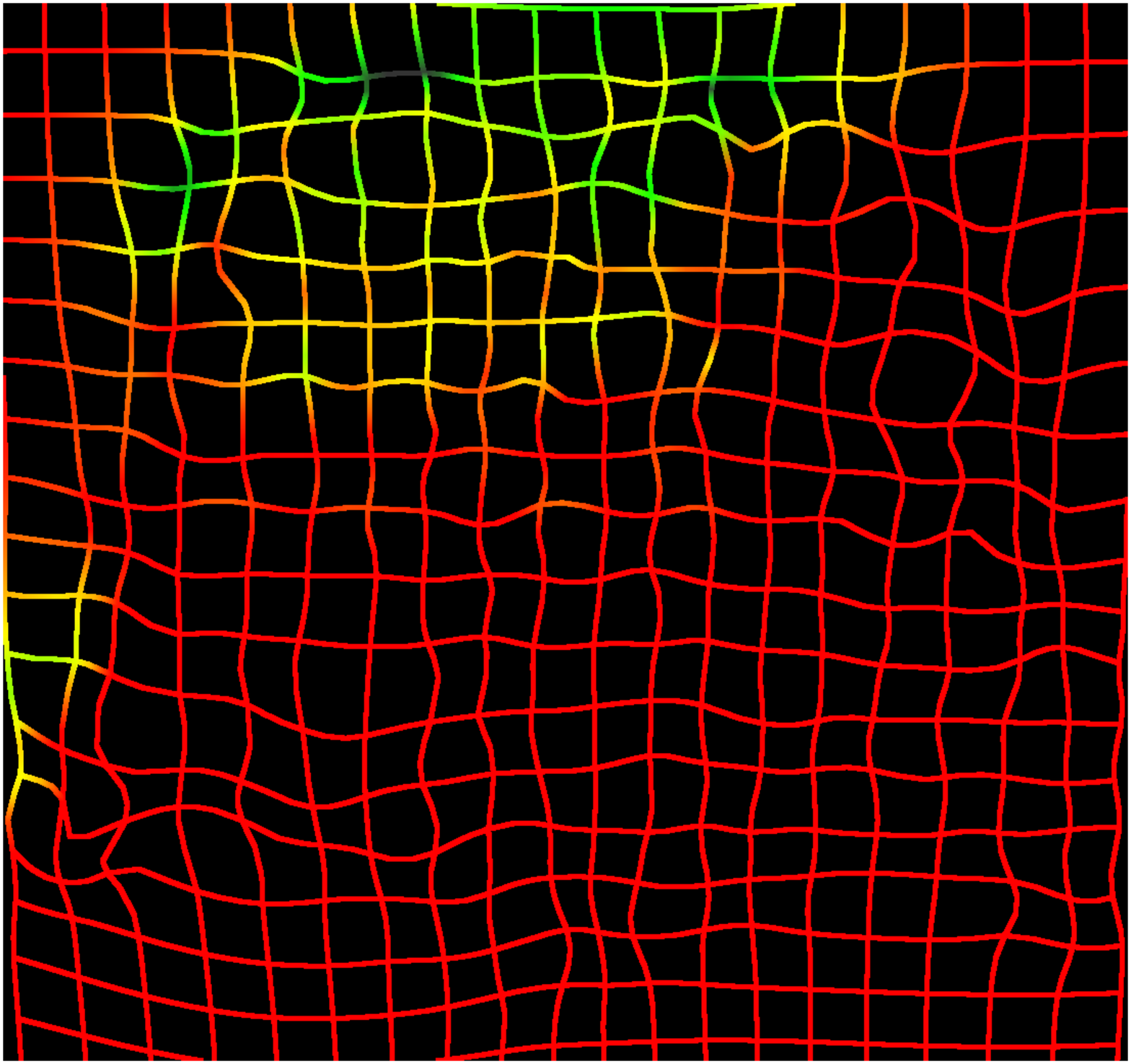}} \hfill
    \subfloat[$\boldsymbol{T}$ for $\mu_2 = 10^{-3}$]{\label{fig:curv}\includegraphics[height=1.5in,width=0.24\textwidth]{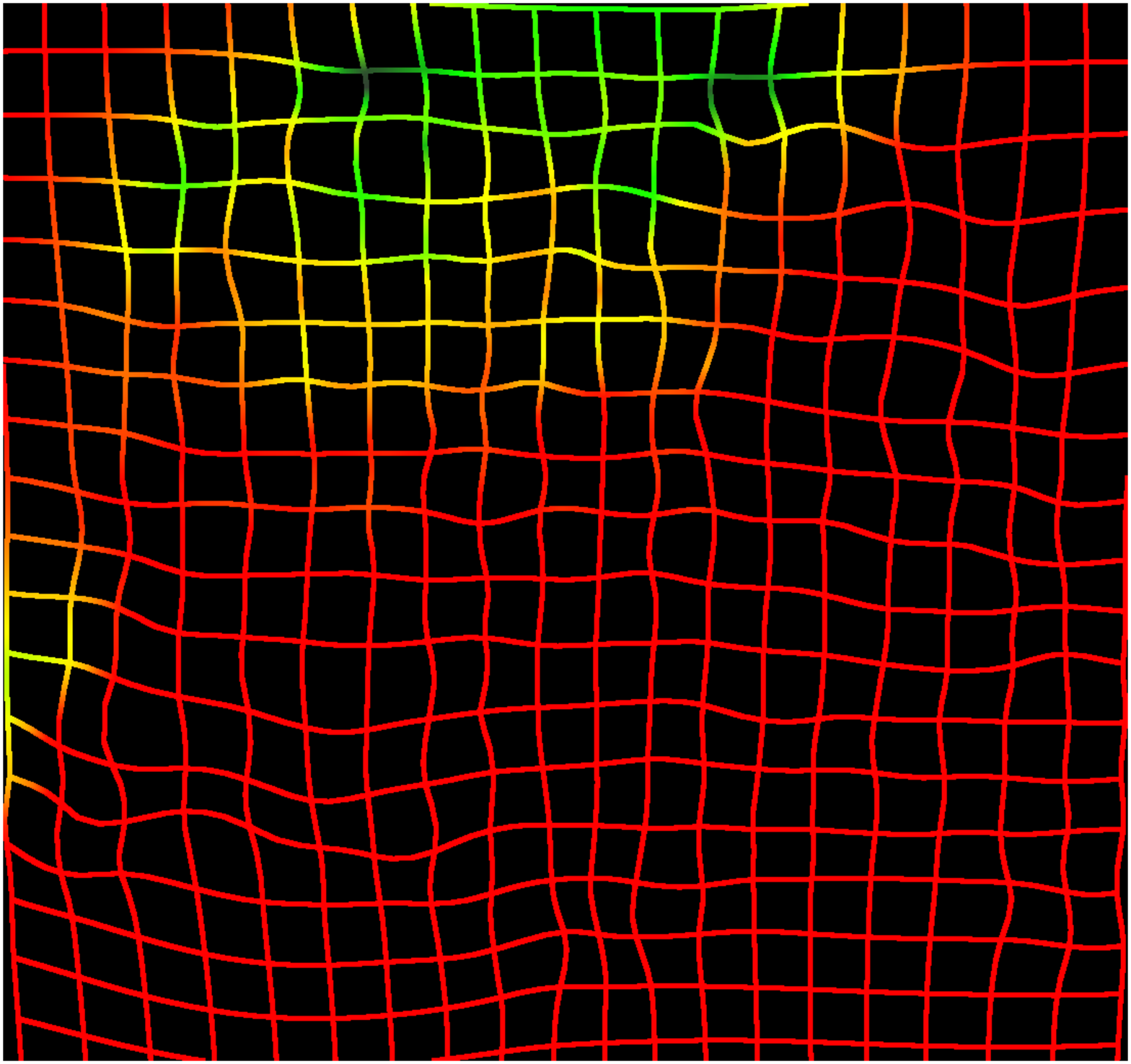}} \hfill
    \subfloat[$\boldsymbol{T}$ for $\mu_2 = 10^{-2}$]{\label{fig:curv_better}\includegraphics[height=1.5in,width=0.24\textwidth]{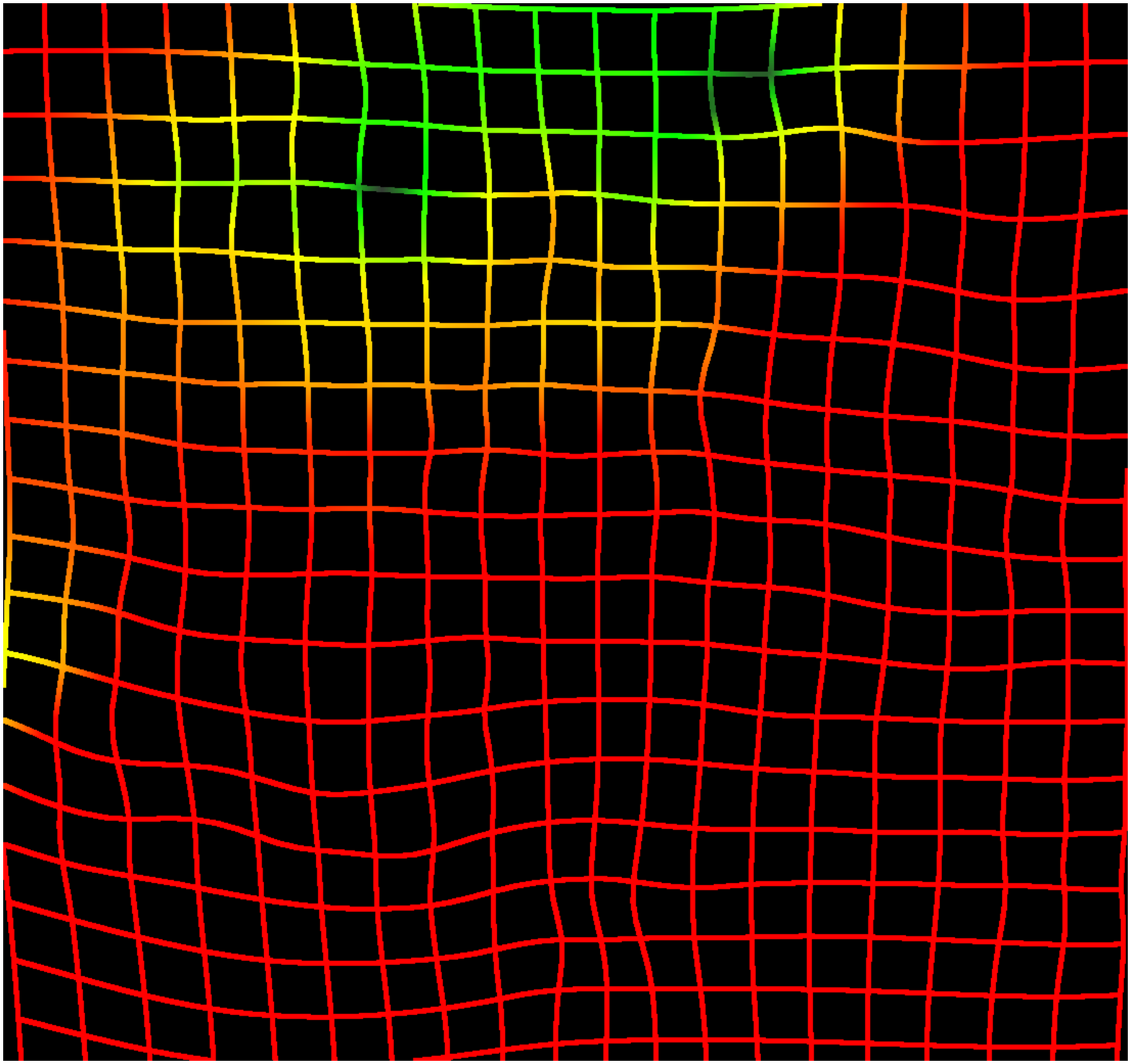}} \hfill
    \caption{Coronal slice showing: (a) the fixed and the moving CT images overlaid; (b) displacement field generated without regularization; (c)--(d) displacement field generated by the curvature regularizer for different values of $\mu_2$. Displacement field shown in black indicates displacement less than 1 mm, green indicates displacement of around 2 mm, yellow indicates displacement of around 5 mm, and red indicates displacement greater than 10 mm, }
    \label{fig:graphs2}
\end{figure*}

The effect of the penalty factor for the curvature regularizer $\mu_2$ on the transformation $\boldsymbol{T}$ is qualitatively shown in Figure~\ref{fig:graphs2} (c-d). The inhaled and exhaled thoracic volumes are registered using a B-spline control-point spacing of $10 \times 10 \times 10$ mm while varying $\mu_n$, and the resulting transform exhibits increased smoothness upon application of the penalty term. Taking the curvature regularizer for example, we see in Figure~\ref{fig:graphs1} (a)-left that for a control point spacing of 10 mm, the MLS is comparable when $\mu_2$ is either $10^{-3}$ or $10^{-2}$. Inspecting the minimum Jacobian determinant in Figure~\ref{fig:graphs1} (b)-left, we see that the minimum Jacobian determinant is negative when $\mu_2$ is $10^{-3}$, which indicates a non-smooth displacement field and is positive when $\mu_2$ is $10^{-2}$, which indicates a smooth displacement field. This makes $\mu_2 = 10^{-2}$ a better choice for the regularization penalty term. This can also be verified by comparing Figures~\ref{fig:graphs2} (c) and~\ref{fig:graphs2} (d). Quantitatively, if the registration accuracy in terms of MLS is the same, the smooth displacement field with a positive minimum Jacobian determinant is preferred. Similar effect is observed for the other four regularizers.

The main advantage of an analytic regularization method is significantly reduced computational time compared to the numerical approach. Computation times incurred by the analytic solutions depend only on the number of tiles defined by the B-spline control-point spacing and not on the number of voxels, drastically reducing the complexity. The composite \textbf{V} matrices can be calculated once for each control point spacing, and reused throughout the optimization. 

\begin{table}[!ht]
    \centering
    \caption{Execution times for the analytic methods and corresponding speedup over numerical approach.}
    \label{tab:my-table2}
    \resizebox{\textwidth}{!}{%
        \begin{tabular}{cccccccc}
            \hline
            Volume Size & Grid Spacing & Analytic (s) & Diffusion & \multicolumn{1}{c}{Curvature} & \multicolumn{1}{c}{Elastic} & \multicolumn{1}{c}{Third-order} & Tot. displacement \\ \hline
            $256 \times 256 \times 94$ & $20 \times 20 \times 20$ & 0.43 & 4.3x & 13.8x & 4.3x & 29.6x & 0.8x \\
            $256 \times 256 \times 94$ & $30 \times 30 \times 30$ & 0.14 & 13.6x & 42.7x & 13.3x & 91.6x & 2.6x \\
            $512 \times 512 \times 128$ & $20 \times 20 \times 20$ & 1.63 & 6.3x & 19.8x & 6.3x & 41.9x & 1.5x \\
            $512 \times 512 \times 128$ & $30 \times 30 \times 30$ & 0.67 & 15.6x & 48.2x & 14.7x & 100.4x & 2.9x \\ \hline
        \end{tabular}%
    }
\end{table}

Table~\ref{tab:my-table2} lists the execution times incurred by the analytic regularization and corresponding speedup compared to numerical solutions based on finite differencing. The benchmarking results reported here were performed using a machine equipped with dual Intel Octo-core Xeon processors clocked at 2.4 GHz and 512 GB of main memory. Examining the third column of the table, note that execution time for a given volume size and grid spacing is the same for any of the analytic regularizers. This is because there is one composite matrix \textbf{V} per partial derivative for a total of 32 such matrices in the generalized framework, all of which are computed irrespective of the regularizer to be used. The regularizer of interest is later selected by setting a non-zero weight to the corresponding term in \eqref{eqn:smooth}. Columns four through eight of the table list speedup over numerical implementations of the various regularizers. We consider single-threaded implementations here. Speedup depends on three factors: volume size, grid spacing, and complexity of the numerical solution. For a volume size of $512 \times 512 \times 128$ voxels with grid spacing of $20 \times 20 \times 20$ mm. The B-spline grid has 7744 tiles; with grid spacing $30 \times 30 \times 30$ mm it has 3177 tiles. In theory, execution time is linear with the number of tiles.
    
Referring back to Table~\ref{tab:my-table2}, third-order and curvature regularizers exhibit higher speedup than the other regularizers due to the high complexity of their numerical solution. For example, for a grid-spacing of 20 mm, the analytic form of the curvature regularizer executes in 450 ms whereas the numerical form requires 6.2 seconds. On the other hand, the analytic form of the total displacement regularizer executes in 450 ms but the numeric form only takes 380 ms because regularization involves simply summing the squared vector-field magnitude at each voxel as per~\eqref{eqn:smooth_td}.

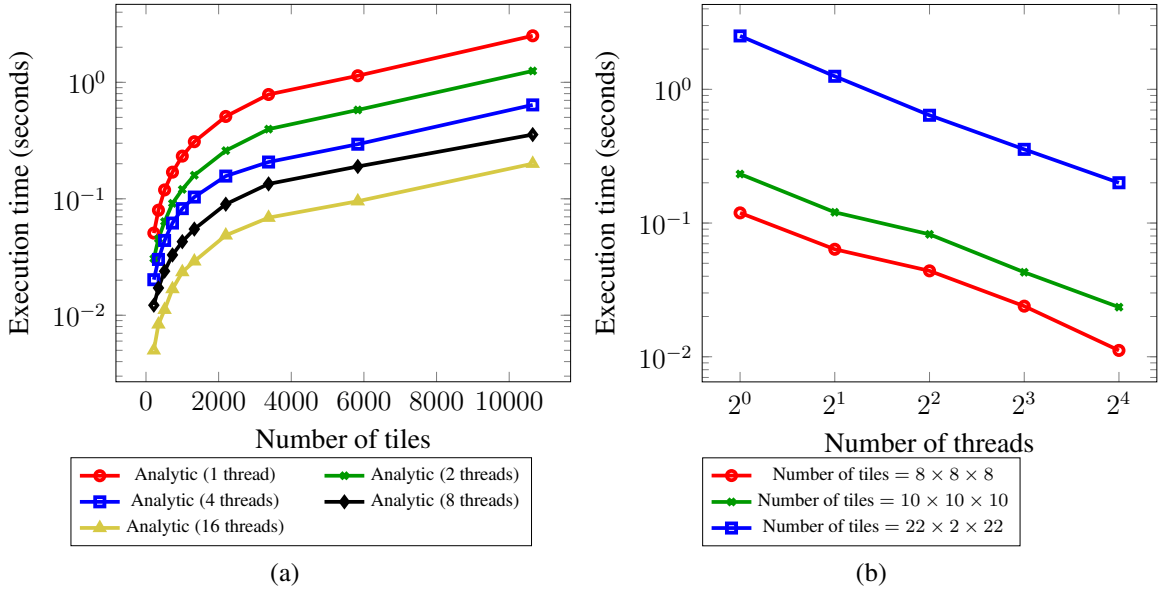
\begin{figure}[!ht]
\centering
\subfloat[]{
\resizebox{0.49\textwidth}{!}{%
    \begin{tikzpicture}
    \begin{semilogyaxis}[
    xlabel={Number of tiles},
    ylabel={Execution time (seconds)},
    every axis plot/.append style={ultra thick},
    legend entries ={Analytic (1 thread),Analytic (2 threads),Analytic (4 threads),Analytic (8 threads),Analytic (16 threads),Analytic (32 threads)},
    legend style={font=\fontsize{7}{7}\selectfont,
	at={(-0.1,-0.2)},
    anchor=north west,
    legend columns=2,
    /tikz/every even column/.append style={column sep=0.5cm}},
    scaled x ticks = false,
    x tick label style={/pgf/number format/.cd,%
        scaled x ticks = false,
        set thousands separator={},
        fixed},]
    \addplot+[red,mark=o] table[line join=round,col sep=comma] {fig5a1.csv};
    \addlegendentry{{\scriptsize Analytic (1 thread)}}
    \addplot+[green!60!black,mark=x] table[line join=round,col sep=comma] {fig5a2.csv};
    \addlegendentry{{\scriptsize Analytic (2 threads)}}
    \addplot+[blue,mark=square] table[line join=round,col sep=comma] {fig5a3.csv};
    \addlegendentry{{\scriptsize Analytic (4 threads)}}
    \addplot+[black,mark=diamond] table[line join=round,col sep=comma] {fig5a4.csv};
    \addlegendentry{{\scriptsize Analytic (8 threads)}}
    \addplot+[yellow!80!black,mark=triangle] table[line join=round,col sep=comma] {fig5a5.csv};
    \addlegendentry{{\scriptsize Analytic (16 threads)}}
    \end{semilogyaxis}
    \end{tikzpicture}
    }}
    \subfloat[]{
    \resizebox{0.49\textwidth}{!}{%
	\begin{tikzpicture}
	\begin{loglogaxis}[
	xlabel={Number of threads},
	ylabel={Execution time (seconds)},
	every axis plot/.append style={ultra thick},
	legend entries ={Number of tiles $=8\times8\times8$, Number of tiles $=10\times10\times10$, Number of tiles $=22\times22\times22$},
	legend style={font=\fontsize{7}{7}\selectfont,
	at={(0.0,-0.2)},
    anchor=north west,
     legend columns=1,
    /tikz/every even column/.append style={column sep=1.0cm}},
	scaled x ticks = false,
	x tick label style={/pgf/number format/.cd,%
		scaled x ticks = false,
		set thousands separator={},
		fixed},
	xtick=data,log basis x=2]
	\addplot+[red,mark=o] table[line join=round,col sep=comma] {fig5b1.csv};
	\addlegendentry{{\scriptsize Number of tiles $=8\times8\times8$}}
	\addplot+[green!60!black,mark=x] table[line join=round,col sep=comma] {fig5b2.csv};
	\addlegendentry{{\scriptsize Number of tiles $=10\times10\times10$}}
	\addplot+[blue,mark=square] table[line join=round,col sep=comma] {fig5b3.csv};
	\addlegendentry{{\scriptsize Number of tiles $=22\times2\times22$}}
	\end{loglogaxis}
	\end{tikzpicture}
	}}

 \caption{Execution times for OpenMP implementations as a function of number of tiles (a) and function of thread count for different numbers of tiles (b).}
    \label{fig:openmp-threads}
\end{figure}
Finally, the analytical formulation developed in the previous section is easily parallelized. Returning to \eqref{eqn:s-final1}, this expression can be calculated for each tile in parallel to calculate partial sums which can then be reduced to a single value to obtain the smoothness penalty for the entire volume. We compare the execution time incurred by the single-threaded implementation against a multi-threaded one using OpenMP. Execution time for a $512 \times 512 \times 512$ image is shown in Figure~\ref{fig:openmp-threads} as as a function of the number of number of tiles in the image as well as the number of threads (from one to sixteen) for the linear-elastic regularizer. Experiments are repeated twenty times to avoid any discrepancy in the timings caused by a cold cache and the Figure~\ref{fig:openmp-threads} shows average execution times. The effect of varying the control-point grid spacing on the execution time of the analytic implementation of the linear-elastic regularizer for different number of threads can be seen in Figure~\ref{fig:openmp-threads} (a). Notice the nearly linear speed-up of about $16x$ when the number of threads is increased to sixteen, which is the number of cores available on the system used to produce these benchmarks. Figure~\ref{fig:openmp-threads} (b) shows variation in execution time as a function of thread count, for specific grid-spacing settings. Execution time decreases when the control-point grid spacing is increased since this reduces the number of tiles in the overall image. Since the underlying implementation of all five regularizer is similar, so no comparison of performance is provided.

\section{Conclusions}\label{sec:conclusions}
We have developed a fast and general framework which supports five unique regularizers to calculate the smoothness penalty using an analytical approach --- specifically, by deriving composite matrix operators that operate on a set of 64 B-spline control points to calculate the regularization penalty within a given region of support. In terms of accuracy, the maximum relative difference between the analytical and numerical solutions was $7.4\%$. Furthermore, the analytic solutions run up to two orders of magnitude faster than finite differencing based numerical solutions. Fast analytical methods such as these provide effective regularization without imposing a computational burden to the deformable image registration pipeline. 
\vspace{-6pt}

\ack{This material is based upon work supported by the National Science Foundation under Grant Nos. 1553436, 1642345 and 1642380 and the National Institutes of Health under NCI R01CA229178.}

\newcommand{\newblock}{}
\bibliographystyle{elsarticle-harv}
\bibliography{BPEX}
\end{document}